\documentclass[12pt,a4paper]{article}
\usepackage{amsmath}
\usepackage{amssymb}
\usepackage{t1enc}
\usepackage[latin1]{inputenc}
\usepackage[english]{babel}
\usepackage{amsfonts}
\usepackage{latexsym}
\usepackage{bm}
\newtheorem{theorem}{Theorem}[section]
\newtheorem{prop}[theorem]{Proposition}
\newtheorem{lemma}[theorem]{Lemma}

\newtheorem{conj}[theorem]{Conjecture}

\newtheorem{rem}[theorem]{Remark}

\begin{document}
\title{Intersecting integer partitions}

\author{Peter Borg\\[5mm]
Department of Mathematics, University of Malta, Malta\\
\texttt{p.borg.02@cantab.net}}
\date{} \maketitle

\begin{abstract}
If $a_1, a_2, ..., a_k$ and $n$ are positive integers 
such that $n = a_1 + a_2 + ... + a_k$, then the sum $a_1 + a_2 + ... + a_k$ is said to be a \emph{partition of $n$} of \emph{length $k$}, and
$a_1, a_2, ..., a_k$ are said to be the \emph{parts} of the
partition. Two partitions that differ only in the order of their parts
are considered to be the same. We say that two
partitions \emph{intersect} if they have at least one common
part. We call a set $A$ of partitions \emph{intersecting} if any
two partitions in $A$ intersect. Let $P_{n,k}$ be the set of all
partitions of $n$ of length $k$. We conjecture that if $2 \leq k
\leq n$, then the size of any intersecting subset of $P_{n,k}$ is
at most the size of $P_{n-1,k-1}$, which is the size of the
intersecting subset of $P_{n,k}$ consisting of those partitions
which have $1$ as a part. The conjecture is trivially true for $n
\leq 2k$, and we prove it for $n \geq 5k^5$. We also generalise
this for subsets of $P_{n,k}$ with the property that any two of their members have at least $t$ common parts.
\end{abstract}

\section{Introduction}\label{Intro}

Unless otherwise stated, we shall use small letters such as $x$ to
denote elements of a set or positive integers or functions, capital letters such as $X$ to denote sets, and calligraphic letters such as $\mathcal{F}$ to denote \emph{families} (i.e.~sets whose elements are sets themselves). 
We call a set $A$ an \emph{$r$-element set} 
if its size $|A|$ is $r$ (i.e.~if it contains exactly $r$ elements). 
The \emph{power set} of a set $X$ (i.e.~the family of all subsets of $X$) is denoted by $2^X$, and the family of all $r$-element subsets of $X$ is denoted by ${X \choose r}$. For any integer $n \geq 1$, the set $\{1, \dots, n\}$ of the first $n$ positive integers is denoted by $[n]$.

In the literature, a sum $a_1 + a_2 + \dots + a_k$ is said to be a \emph{partition of $n$} of \emph{length $k$} if $a_1, a_2, \dots, a_k$ and $n$ are positive integers such that $n = a_1 + a_2 + \dots + a_k$. 
If $a_1 + a_2 + \dots + a_k$ is a partition, then $a_1, a_2, ..., a_k$ are said to be its \emph{parts}. Two partitions that differ only in the order of their parts are considered to be the same. 
Thus, we can refine the definition of a partition as follows. 
We call a tuple $(a_1, \dots, a_k)$ a \emph{partition of $n$} of \emph{length $k$} if $a_1, \dots, a_k$ and $n$ are positive integers such that $n = \sum_{i=1}^k a_i$ and $a_1 \leq \dots \leq a_k$. We will be using the latter definition throughout the rest of the paper.

For any tuple ${\bf a} = (a_1, \dots, a_k)$ and any $i \in [k]$, we call $a_i$ the \emph{$i$'th entry of ${\bf a}$}, and if $\bf a$ is a partition, then we also call $a_i$ a \emph{part of $\bf a$}.

For any $n$, let $P_n$ be the set of all partitions of $n$, and for any $k$, let $P_{n,k}$ be the set of all partitions of $n$ of length $k$. So $P_{n,k}$ is non-empty if and only if $1 \leq k \leq n$. Moreover, $P_n = \bigcup_{i=1}^n P_{n,i}$. Let $p_n = |P_{n}|$ and $p_{n,k} = |P_{n,k}|$. To the best of the author's knowledge, no closed-form expression is known for $p_n$ and $p_{n,k}$; for more about these values, we refer the reader to \cite{AE}. 

For any set $A$ of integer partitions and any set $T$ of positive integers, let $A(T)$ denote the set of all members of $A$ which have each integer in $T$ as a part; so $A(T) = \{{\bf a} \in A \colon T \mbox{ is a subset of the set of parts of } {\bf a}\}$. For an integer $a$, we may abbreviate the notation $A(\{a\})$ to $A(a)$. Thus, we have
\[P_{n,k}(1) = \{(a_1, \dots, a_k) \in P_{n,k} \colon a_1 = 1\} \quad \mbox{and} \quad P_n(1) = \bigcup_{i=1}^n P_{n,i}(1).\] 
Note that $|P_n(1)| = p_{n-1}$ and $|P_{n,k}(1)| = p_{n-1,k-1}$. 

We say that two partitions \emph{intersect} if they have at least one common part; in other words, a partition $(a_1, \dots, a_r)$ intersects a partition $(b_1, \dots, b_s)$ if $a_i = b_j$ for some $i \in [r]$ and $j \in [s]$. We call a set $A$ of partitions \emph{intersecting} if any two partitions in $A$ intersect. Thus, for any set $A$ of partitions, $A(1)$ is intersecting.

We suggest the following two conjectures.


\begin{conj}[Weak Form] \label{wf} $P_n(1)$ is an intersecting subset of $P_n$ of maximum size.
\end{conj}

\begin{conj}[Strong Form] \label{sf} For $2 \leq k \leq n$, $P_{n,k}(1)$ is an intersecting subset of $P_{n,k}$ of maximum size.
\end{conj}

\begin{prop}\label{prop1} If Conjecture~\ref{sf} is true, then Conjecture~\ref{wf} is true, and for any $n \geq 3$, $P_n(1)$ is the unique intersecting subset of $P_n$ of maximum size.
\end{prop} 
\textbf{Proof.} The result is trivial for $n \leq 2$, so suppose $n \geq 3$. Let $A$ be an intersecting subset of $P_{n}$. For each $k \in [n]$, let $A_k = A \cap P_{n,k}$. So $A_1, \dots, A_n$ are intersecting, and $|A| = \sum_{k = 1}^n |A_k|$. $(n)$ is the only partition in $P_{n,1}$. No partition in $P_n$ other than $(n)$ intersects $(n)$. Thus, if $(n) \in A$, then $A = \{(n)\}$ and hence $|A| = 1 < |P_n(1)|$. Now suppose $(n) \notin A$. Then $A_1 = \emptyset$. Suppose Conjecture~\ref{sf} is true. Then $|A_k| \leq |P_{n,k}(1)|$ for any $k \in [n]$. So we have $|A| = \sum_{k = 2}^n |A_k| \leq \sum_{k = 2}^n |P_{n,k}(1)| = |P_n(1)|$. $P_{n,n}$ has only one partition $\bf a$, and ${\bf a} = (1, \dots, 1)$. If ${\bf a} \in A$, then, since $A$ is intersecting, $A \subseteq P_{n}(1)$. If ${\bf a} \notin A$, then $A_n = \emptyset$ and hence $|A| = \sum_{k = 2}^{n-1} |A_k| \leq \sum_{k = 2}^{n-1} |P_{n,k}(1)| < \sum_{k = 2}^{n} |P_{n,k}(1)| = |P_n(1)|$.~\hfill{$\Box$}\\


We will prove that if $k \geq 3$ and $n$ is sufficiently larger than $k$, then $P_{n,k}(1)$ is the unique intersecting subset of $P_{n,k}$ of maximum size. We suspect that this holds for any $n \geq k \geq 4$. This is not the case for $k = 2$ and $n \geq 4$ (see below), and for $k = 3$ and some values of $n$ (it is easy to check this for $6 \leq n \leq 10$). For example, $\{(1,2,7), (1,3,6), (1,4,5), (2,3,5)\}$ is an intersecting subset of $P_{10,3}$ of size $|P_{10,3}(1)|$. 

Conjecture~\ref{sf} is true for $n \leq 2k$.

\begin{prop}\label{prop2} Let $2 \leq k \leq n \leq 2k$. Then $P_{n,k}(1)$ is an intersecting subset of $P_{n,k}$ of maximum size, and uniquely so unless $2 \leq k \leq 3$ and $n = 2k$.
\end{prop}
\textbf{Proof.} The result is trivial for $2 \leq k \leq 3$, so consider $k \geq 4$. Suppose $n < 2k$. Then every partition of $n$ of length $k$ has $1$ as a part (because the sum of $k$ integers that are all greater than $1$ is at least $2k$). Now suppose $n=2k$. Let ${\bf a}_1$ be the partition in $P_{n,k}$ whose $k$ entries are all $2$. Then ${\bf a}_1$ is the only partition in $P_{n,k}$ that does not have $1$ as a part. Let ${\bf a}_2$ be the partition $(a_1, \dots, a_k)$ in $P_{n,k}(1)$ with $a_1 = \dots = a_{k-1} = 1$ and $a_k = n - k + 1$. So ${\bf a}_1$ does not intersect ${\bf a}_2$. Suppose $A$ is an intersecting subset of $P_{n,k}$ that contains ${\bf a}_1$. Then ${\bf a}_2 \notin A$ and hence $|A| \leq |P_{n,k}(1)|$. Let ${\bf a}_3 = (b_1, \dots, b_k)$ with $b_1 = \dots = b_{k-2} = 1$, $b_{k-1} = 3$ and $b_k = n-k-1$. Then ${\bf a}_3 \in P_{n,k}(1)$ and ${\bf a}_3 \notin A$. So $|A| < |P_{n,k}(1)|$.~\hfill{$\Box$}\\

Using ideas from \cite{Borg}, we will prove that Conjecture~\ref{sf} is also true for $n \geq 5k^5$. 

\begin{theorem}\label{main} For $k \geq 3$ and $n \geq 5k^5$, $P_{n,k}(1)$ is the unique intersecting subset of $P_{n,k}$ of maximum size.
\end{theorem}
We will actually prove two generalisations (Theorems~\ref{main2} and \ref{main3}) of this result in Section~\ref{proofsection}. For $k=2$, $P_{n,k}(1)$ is a largest intersecting subset of $P_{n,k}$, but not uniquely so if $n \geq 4$. 
Indeed, each partition ${\bf a} = (a_1, a_2)$ in $P_{n,2}$ must have $a_2 = n-a_1$, and hence no other partition in $P_{n,2}$ intersects $\bf a$; so an intersecting subset of $P_{n,2}$ cannot have more than one member.

Theorem~\ref{main} is an analogue of the classical Erd\H os-Ko-Rado (EKR) Theorem \cite{EKR}, which inspired many results in extremal set theory (see \cite{Borg7,DF,F}). With a slight abuse of terminology, we say that a family $\mathcal{A}$ of sets is \emph{intersecting} if any two sets in $\mathcal{A}$ intersect (i.e.~$A \cap B \neq \emptyset$ for any $A, B \in \mathcal{A}$). The EKR Theorem says that if $r \leq n/2$ and $\mathcal{A}$ is an intersecting subfamily of ${[n] \choose r}$, then $|\mathcal{A}| \leq {n-1 \choose r-1}$, and equality holds if $\mathcal{A} = \{A \in {[n] \choose r} \colon 1 \in A\}$. 


\begin{rem}\label{rem1} \emph{The above conjectures and results can be re-phrased in terms of intersecting subfamilies of a family. For any integer partition ${\bf a} = (a_1, \dots, a_k)$, let $S_{\bf a} = \{(a,i) \colon i \in [k], \, |\{j \in [k] \colon a_j = a\}| \geq i\}$; thus, $(a,1), \dots, (a,r) \in S_{\bf a}$ if and only if $r$ of the entries of $\bf a$ are $a$. For example, $S_{(2,2,5,5,5,7)} = \{(2,1), (2,2), (5,1), (5,2), (5,3), (7,1)\}$. Let $\mathcal{P}_n = \{S_{\bf a} \colon {\bf a} \in P_{n}\}$ and $\mathcal{P}_{n,k} = \{S_{\bf a} \colon {\bf a} \in P_{n,k}\}$. Let $f : P_{n} \rightarrow \mathcal{P}_{n}$ such that $f({\bf a}) = S_{\bf a}$ for each ${\bf a} \in P_{n,k}$. Clearly, $f$ is a bijection. So $|\mathcal{P}_n| = |P_n|$ and $|\mathcal{P}_{n,k}| = |P_{n,k}|$. Note that two integer partitions ${\bf a}$ and $\bf b$ intersect if and only if $S_{\bf a} \cap S_{\bf b} \neq \emptyset$. Thus, a subset $A$ of $P_{n,k}$ is a largest intersecting subset if and only if $\{S_{\bf a} \colon {\bf a} \in A\}$ is a largest intersecting subfamily of $\mathcal{P}_{n,k}$.} 
\end{rem}

\section{$t$-intersecting integer partitions}\label{tint}

A family $\mathcal{A}$ of sets is said to be \emph{$t$-intersecting} if $|A \cap B| \geq t$ for any $A, B \in \mathcal{A}$. A $t$-intersecting family $\mathcal{A}$ is said to be \emph{non-trivial} if $|\bigcap_{A \in \mathcal{A}} A| < t$ (i.e.~the number of elements common to all the sets in $\mathcal{A}$ is less than $t$). 
Note that an intersecting family is a $1$-intersecting family. In addition to the EKR Theorem (see Section~\ref{Intro}), it was also proved in \cite{EKR} that if $n$ is sufficiently larger than $r$, then the size of any $t$-intersecting subfamily of ${[n] \choose r}$ is at most ${n-t \choose r-t}$, and hence $\{A \in {[n] \choose r} \colon [t] \subset A\}$ is a largest $t$-intersecting subfamily of ${[n] \choose r}$. The complete solution for any $n$, $r$ and $t$ is given in \cite{AK1}; it turns out that $\{A \in {[n] \choose r} \colon [t] \subset A\}$ is a largest $t$-intersecting subfamily of ${[n] \choose r}$ if and only if $n \geq (r-t+1)(t+1)$.

We introduce two generalisations of the definition of an intersecting set of integer partitions. 

Let ${\bf a} = (a_1, \dots, a_r)$ and ${\bf b} = (b_1, \dots, b_s)$ be two integer partitions. We say that $\bf a$ and $\bf b$ \emph{$t$-intersect} if they have $t$ common parts (not necessarily distinct); more precisely, $\bf a$ $t$-intersects $\bf b$ if there are $t$ distinct integers $i_1, \dots, i_t$ in $[r]$ and $t$ distinct integers $j_1, \dots, j_t$ in $[s]$ such that $a_{i_p} = b_{j_p}$ for each $p \in [t]$. We say that $\bf a$ and $\bf b$ \emph{$t$-intersect properly} if they have $t$ distinct common parts; in other words, $\bf a$ $t$-intersects $\bf b$ properly if $|\{a_i \colon i \in [r]\} \cap \{b_j \colon j \in [s]| \geq t$. Note that if $\bf a$ and $\bf b$ $t$-intersect properly, then $\bf a$ and $\bf b$ $t$-intersect.

Let $A$ be a set of integer partitions. With a slight abuse of terminology, we say that $A$ is \emph{$t$-intersecting} if for any ${\bf a}, {\bf b} \in A$, $\bf a$ and $\bf b$ $t$-intersect. We say that $A$ is \emph{properly $t$-intersecting} if for any ${\bf a}, {\bf b} \in A$, $\bf a$ and $\bf b$ $t$-intersect properly. Note that an intersecting set of integer partitions is $1$-intersecting and properly $1$-intersecting. 

We suggest generalisations of Conjectures~\ref{wf} and \ref{sf} along the lines of the general definitions above.

For any set $A$ of integer partitions, let $A \langle t \rangle$ denote the set of all partitions in $A$ whose first $t$ entries are $1$. Thus, for $1 \leq t \leq k \leq n$, 
\[P_{n,k} \langle t \rangle = \{(a_1, \dots, a_k) \in P_{n,k} \colon a_1 = \dots = a_t = 1\}\} \quad \mbox{and} \quad P_n\langle t \rangle = \bigcup_{i=t}^n P_{n,i}\langle t \rangle.\] 
Note that $|P_n\langle t \rangle| = |P_{n-t}|$ and $|P_{n,k}\langle t \rangle| = |P_{n-t,k-t}|$. Also note that $A(1) = A \langle 1 \rangle$.

We first present and discuss our conjectures for $t$-intersecting partitions.

\begin{conj}\label{tconj1a} $P_{n}\langle t \rangle$ is a $t$-intersecting subset of $P_{n}$ of maximum size.
\end{conj}

\begin{conj}\label{tconj1b} For $t + 1 \leq k \leq n$, $P_{n,k}\langle t \rangle$ is a $t$-intersecting subset of $P_{n,k}$ of maximum size.
\end{conj}
Note that if $t = k < n$, then $P_{n,k}\langle t \rangle = \emptyset$, $P_{n,k} \neq \emptyset$, and the $t$-intersecting subsets of $P_{n,k}$ are the $1$-element subsets. If $k < t$, then $P_{n,k}$ has no non-empty $t$-intersecting subsets. If Conjecture~\ref{tconj1b} is true, then, by an argument similar to that of Proposition~\ref{prop1}, Conjecture~\ref{tconj1a} is true.

\begin{prop} Conjecture~\ref{tconj1b} is true for $n \leq 2k-t+1$.
\end{prop}
\textbf{Proof.} By Proposition~\ref{prop2}, we may assume that $t \geq 2$. Suppose $n \leq 2k - t+1$. For any ${\bf c} = (c_1, \dots, c_k) \in P_{n,k}$, let $L_{\bf c} = \{i \in [k] \colon c_i = 1\}$ and let $l_{\bf c} = |L_{\bf c}|$. 

Let ${\bf c} = (c_1, \dots, c_k) \in P_{n,k}$. We have $2k - t + 1 \geq n = \sum_{i \in L_{\bf c}}c_i + \sum_{j \in [k] \backslash L_{\bf c}} c_j \geq \sum_{i \in L_{\bf c}} 1 + \sum_{j \in [k] \backslash L_{\bf c}} 2 = l_{\bf c} + 2(k - l_{\bf c}) = 2k - l_{\bf c}$. Thus, $l_{\bf c} \geq t-1$, and equality holds only if $n = 2k-t+1$ and $c_j = 2$ for each $j \in [k] \backslash L_{\bf c}$. Since $c_1 \leq \dots \leq c_k$, $L_{\bf c} = [l_{\bf c}]$.

Let $A$ be a $t$-intersecting subset of $P_{n,k}$. If $l_{\bf a} \geq t$ for each ${\bf a} \in A$, then $A \subseteq P_{n,k}\langle t \rangle$. Suppose $l_{\bf a} = t-1$ for some ${\bf a} = (a_1, \dots, a_k) \in A$. Then, by the above, we have $n = 2k-t+1$, $a_i = 1$ for each $i \in [t-1]$, $a_j = 2$ for each $j \in [k] \backslash [t-1]$, and $P_{n,k} = P_{n,k}\langle t \rangle \cup \{{\bf a}\}$. Let ${\bf b}$ be the partition $(b_1, \dots, b_k)$ in $P_{n,k}\langle t \rangle$ with $b_k = n - k + 1 = k - t + 2$ and $b_i = 1$ for each $i \in [k-1]$. So $\bf a$ and $\bf b$ do not $t$-intersect, and hence ${\bf b} \notin A$. So $|A| \leq |P_{n,k}| - 1 = |P_{n,k}\langle t \rangle|$.~\hfill{$\Box$}\\

The following generalisation of Theorem~\ref{main} tells us that Conjecture~\ref{tconj1b} is also true for $n \geq {3k-2t-1 \choose t+1}k^3$. Its proof is given in the next section.

\begin{theorem}\label{main2} For $k \geq t + 2$ and $n \geq {3k-2t-1 \choose t+1}k^3$, $P_{n,k}\langle t \rangle$ is the unique $t$-intersecting subset of $P_{n,k}$ of maximum size.
\end{theorem}
Conjecture~\ref{tconj1b} is also true for $k = t+1$. Indeed, if two partitions of $n$ of length $t+1$ have $t$ common parts $a_1, \dots, a_t$ (not necessarily distinct), then the remaining part of each partition is $n - (a_1 + \dots + a_t)$, and hence the partitions are the same. Thus, the $t$-intersecting subsets of $P_{n,t+1}$ are the $1$-element subsets. So $P_{n,t+1}\langle t \rangle$ is a largest $t$-intersecting subset of $P_{n,k}$, but not uniquely so if $n \geq t+3$ (because in this case, at least there is also $\{(1, \dots, 1, 2, n-t-1)\}$). 

We now present and discuss our conjectures for properly $t$-intersecting partitions.

\begin{conj}\label{tconj2a} $P_{n}([t])$ is a properly $t$-intersecting subset of $P_{n}$ of maximum size.
\end{conj}

\begin{conj}\label{tconj2b} For $t + 1 \leq k \leq n$, $P_{n,k}([t])$ is a properly $t$-intersecting subset of $P_{n,k}$ of maximum size.
\end{conj}
Conjecture~\ref{tconj2a} is trivial for $n < t(t+1)/2$, and Conjecture~\ref{tconj2b} is trivial for $n < t(t-1)/2 + k$. Indeed, each member of a properly $t$-intersecting set of partitions must have at least $t$ distinct parts; thus, $P_{n,k}$ has no non-empty $t$-intersecting subsets if $n < \sum_{i=1}^t i = t(t+1)/2$, and $P_{n,k}$ has no non-empty $t$-intersecting subsets if $n < t(t+1)/2 + k-t = t(t-1)/2 + k$. The reason why we need $k \geq t+1$ in Conjecture~\ref{tconj2b} is similar to that for Conjecture~\ref{tconj1b}. If Conjecture~\ref{tconj2b} is true, then, by an argument similar to that of Proposition~\ref{prop1}, Conjecture~\ref{tconj2a} is true.

The following generalisation of Theorem~\ref{main} tells us that Conjecture~\ref{tconj2b} is true for $n \geq {3k-2t-1 \choose t+1}k^3$. Its proof is given in the next section.

\begin{theorem}\label{main3} For $k \geq t + 2$ and $n \geq {3k-2t-1 \choose t+1}k^3$, $P_{n,k}([t])$ is the unique properly $t$-intersecting subset of $P_{n,k}$ of maximum size.
\end{theorem}
Similarly to Conjecture~\ref{tconj1b}, Conjecture~\ref{tconj2b} is also true for $k = t+1$, but $P_{n,t+1}([t])$ is not the unique largest properly $t$-intersecting subset of $P_{n,t+1}$ for $n \geq t(t+1)/2 + 2$ when $t \geq 2$ (indeed, $\{(1, 1, 2, \dots, t-1, t+1)\}$ is another one).

We now proceed by proving Theorems~\ref{main2} and \ref{main3}.

\section{Proofs of Theorems~\ref{main2} and \ref{main3}}\label{proofsection}

\begin{lemma} \label{lemma1} Let $1 \leq k \leq m \leq n$. Then $p_{m,k} \leq p_{n,k}$. Moreover, if $n > m$, $n \geq k+2$ and $k \geq 3$, then $p_{m,k} < p_{n,k}$.
\end{lemma}
\textbf{Proof.} The case $k = 1$ is trivial. Suppose $k \geq 2$.
Let $f \colon P_{m,k} \rightarrow P_{n,k}$ be
the function that maps any partition $(a_1, \dots, a_k)$ in
$P_{m,k}$ to the partition $(b_1, \dots, b_k)$ in
$P_{n,k}$ with $b_k = a_k + n-m$ and $b_i = a_i$ for
any $i \in [k-1]$. Clearly, $f$ is one-to-one, and hence the size of its domain $P_{m,k}$ is at most the size of its co-domain $P_{n,k}$. So $p_{m,k} \leq p_{n,k}$.

Suppose $n > m$, $n \geq k+2$ and $k \geq 3$. Let ${\bf c} = (c_1, \dots, c_k)$ with $c_i = 1$ for each $i \in [k] \backslash \{k-2, k-1, k\}$, $c_{k-2} = 1$ and $c_{k-1} = c_k = (n-k+2)/2$ if $n-k$ is even, and $c_{k-2} = 2$ and $c_{k-1} = c_k = (n-k+1)/2$ if $n-k$ is odd. 
%
%
So ${\bf c} \in P_{n,k}$. Since $m < n$, $f$ maps each partition $(a_1, \dots, a_k)$ in $P_{n,k}$ to a partition $(b_1, \dots, b_k)$ with $b_{k-1} < b_k$. So ${\bf c}$ is not in the range of $f$. So $f$ is not onto, and hence the size of its domain $P_{m,k}$ is less than the size of its co-domain $P_{n,k}$. So $p_{m,k} < p_{n,k}$.~\hfill{$\Box$}

\begin{lemma}\label{lemma2} Let $c \geq 1$ and $k \geq 2$. For any $n
\geq ck^3$, \[p_{n,k} > cp_{n,k-1}.\]
\end{lemma}
\textbf{Proof.} The result is trivial for $k = 2$, so we assume
$k \geq 3$. For each $i \in [ck^2]$, let $F_i = \{(i, a_1, ...,
a_{k-2}, a_{k-1} - i) \colon (a_1, ...,a_{k-1}) \in P_{n,k-1}\}$.
Let $F = \bigcup_{i=1}^n F_i$.

For any $k$-tuple ${\bf x} = (x_1, ..., x_k)$ of positive integers,
let ${\bf x}^{\rightarrow}$ be the $k$-tuple obtained by putting
the entries of ${\bf x}$ in increasing order; that is, ${\bf x}^{\rightarrow}$ is the $k$-tuple $(x_1', \dots, x_k')$ such that $x_1' \leq \dots \leq x_k'$ and, for each $j \in [k]$, $|\{i \in [k] \colon x_i = x_j\}| = |\{i \in [k] \colon x_i' = x_j\}|$. 

Let $\bf a$ be a partition $(a_1, ..., a_{k-1})$ in $P_{n,k-1}$. Since $a_1 \leq \dots \leq a_{k-1}$ and $a_1 + \dots + a_{k-1} = n$, we have $a_{k-1} \geq \frac{n}{k-1}$, and hence, since $n \geq ck^3$, $a_{k-1}
> ck^2$. So $a_{k-1} - i \geq 1$ for each $i \in [ck^2]$, meaning
that the entries of each member of $F$ are positive integers that
add up to $n$. Therefore,
\begin{equation} {\bf x}^{\rightarrow} \in P_{n,k} \mbox{ for each
} {\bf x} \in F. \label{FG}
\end{equation}

Let $G = \{{\bf y} \in P_{n,k} \colon {\bf y} = {\bf
x}^{\rightarrow} \mbox{ for some } {\bf x} \in F\}$. For each
${\bf y} \in G$, let $F_{{\bf y}} = \{{\bf x} \in F \colon {\bf
x}^{\rightarrow} = {\bf y}\}$. By (\ref{FG}), $F \subseteq
\bigcup_{{\bf y} \in G} F_{{\bf y}}$.

Let ${\bf y}$ be a partition $(y_1, ..., y_k)$ in $G$. Clearly, each member of $F_{\bf y}$ is in one of $F_{y_1}, \dots, F_{y_k}$; that is, $F_{\bf y} \subseteq \bigcup_{i=1}^k F_{y_i}$. So $F_{\bf y} = \bigcup_{i=1}^k \left(F_{\bf y} \cap F_{y_i} \right)$. Let $i \in [k]$ such that $F_{\bf y} \cap F_{y_i} \neq \emptyset$. Let ${\bf x}$ be a tuple $(x_1, \dots, x_k)$ in $F_{\bf y} \cap F_{y_i}$. By definition, $x_1 = y_i$ and $x_2 \leq \dots \leq x_{k-1}$. Thus, since ${\bf y} = {\bf x}^{\rightarrow}$ and $y_1 \leq \dots \leq y_k$, $\bf x$ is one of $(y_i, y_2, \dots, y_{i-1}, y_{i+1}, \dots,  y_k)$, $(y_i, y_2, \dots,  y_{i-1}, y_{i+1}, \dots, y_{k-2}, y_k, y_{k-1})$, ..., $(y_i, y_3,  \dots,  y_{i-1}, y_{i+1}, \dots, y_k, y_2)$ (i.e.~the $k-1$ $k$-tuples satisfying the following:~the first entry is $y_i$, the $k$'th entry is $y_j$ for some $j \in [k] \backslash \{i\}$, and the middle $k-2$ entries form the $(k-2)$-tuple obtained by deleting the $i$'th and $j$'th entry from $\bf y$). So $|F_{\bf y} \cap F_{y_i}| \leq k-1$.

Therefore, we have
\begin{align} |F| &= \left|\bigcup_{{\bf y} \in G} F_{{\bf y}} \right| \leq \sum_{{\bf y} \in G} |F_{{\bf y}}| = \sum_{{\bf y} \in G} \sum_{i=1}^k |F_{{\bf y}} \cap F_{y_i}| \leq \sum_{{\bf y} \in G} \sum_{i=1}^k (k-1) = |G|k(k-1) < |P_{n,k}|k^2 \nonumber
\end{align}
and hence $p_{n,k} > \frac{|F|}{k^2}$. Now $F_1, \dots, F_{ck^2}$
are disjoint sets, each of size $p_{n,k-1}$. So $|F| =
ck^2p_{n,k-1}$ and hence $p_{n,k} > cp_{n,k-1}$.~\hfill{$\Box$}
\\

The last lemma we need before proving Theorems~\ref{main2} and \ref{main3} emerges from \cite{EKR}. 

\begin{lemma}\label{lemma3} Let $\mathcal{A}$ be a non-trivial $t$-intersecting family such that $|A| \leq r$ for any $A \in \mathcal{A}$. Then there exists a set $J$ of size at most $3r-2t-1$ such that $|A \cap J| \geq t+1$ for any $A \in \mathcal{A}$.
\end{lemma}
\textbf{Proof.} If $\mathcal{A}$ is $(t+1)$-intersecting, then we just take $J$ to be an arbitrary set in $\mathcal{A}$. So suppose $\mathcal{A}$ is not $(t+1)$-intersecting. Then there exist $A_1, A_2 \in \mathcal{A}$ such that $|A_1 \cap A_2| = t$. Thus, since $\mathcal{A}$ is a non-trivial $t$-intersecting family, there exists $A_3 \in \mathcal{A}$ such that $A_1 \cap A_2 \nsubseteq A_3$, and hence $|A_1 \cap A_2 \cap A_3| \leq t-1$. Take $J$ to be $A_1 \cup A_2 \cup A_3$. So $|A \cap J| \geq t$ for all $A \in \mathcal{A}$. Suppose there exists $A^* \in \mathcal{A}$ such that $|A^* \cap J| = t$. Then $t \geq |A^* \cap (A_1 \cup A_2)| = |A^*
\cap A_1| + |A^* \cap A_2| - |A^* \cap A_1 \cap A_2| \geq 2t - |A^* \cap A_1 \cap A_2|$, and hence $t \leq |A^* \cap A_1 \cap A_2|$. Since $|A^* \cap A_1 \cap A_2| \leq |A^* \cap J| = t$, we actually have $|A^* \cap A_1 \cap A_2| = |A^* \cap J|$, and hence $A^* \cap J = A^* \cap A_1 \cap A_2$ (as $A_1 \cap A_2 \subset J$). So we have $t \leq |A^* \cap A_3| = |A^* \cap (A_3 \cap J)| = |(A^* \cap J) \cap A_3| = |(A^* \cap A_1 \cap A_2) \cap A_3| \leq |A_1 \cap A_2 \cap A_3|$, which contradicts $|A_1 \cap A_2 \cap A_3| \leq t-1$. So $|A \cap J| \geq t+1$ for all $A
\in \mathcal{A}$. Now $|J| = |A_1 \cup A_2| + |A_3| - |A_3 \cap (A_1 \cup A_2)|$. Since $|A_1 \cup A_2| \leq 2r - |A_1 \cap A_2| = 2r-t$ and $|A_3 \cap (A_1 \cup A_2)| = |A_3 \cap A_1| + |A_3 \cap A_2| - |A_3 \cap A_2 \cap A_1| \geq 2t - |A_1 \cap A_2 \cap A_3| \geq 2t - (t-1) = t+1$, it follows that $|J| \leq (2r-t) + r - (t+1) = 3r-2t-1$.~\hfill{$\Box$}
\\
\\
\textbf{Proof of Theorem~\ref{main2}.} 
Let $k \geq t+2$ and $n \geq {3k-2t-1 \choose t+1}k^3$. Let $A$ be an intersecting subset of $P_{n,k}$ such that $A \neq P_{n,k}\langle t \rangle$. We prove the result by showing that $|A| < |P_{n,k}\langle t \rangle|$.

For each ${\bf a} \in P_{n,k}$, let $S_{\bf a}$ be as in Remark~\ref{rem1}. Define $\mathcal{P}_{n,k}$ and $f$ also as in Remark~\ref{rem1}. Let $\mathcal{A} = \{f({\bf a}) \colon {\bf a} \in A\}$. Clearly, $|\mathcal{A}| = |A|$ (since $f$ is a bijection) and $|X| = k$ for each $X \in \mathcal{A}$. Note that two integer partitions ${\bf a}$ and $\bf b$ in $P_{n,k}$ $t$-intersect if and only if $|S_{\bf a} \cap S_{\bf b}| \geq t$. Thus, since $A$ is a $t$-intersecting set, $\mathcal{A}$ is a $t$-intersecting family. 

Suppose the sets in $\mathcal{A}$ have $t$ common elements $(c_1,i_1), \dots, (c_t,i_t)$. Then $t$ of the entries of each member of $A$ are $c_1, \dots, c_t$. So $|A| \leq p_{n-s,k-t}$, where $s = c_1 + \dots + c_t$. Now $s \geq t$. If $s > t$, then, by Lemma~\ref{lemma1}, $|A| < p_{n-t,k-t} = |P_{n,k}\langle t \rangle|$. Suppose $s = t$. Then $c_1 = \dots = c_t = 1$. So $A \subset P_{n,k}\langle t \rangle$. Since $A \neq P_{n,k}\langle t \rangle$, $|A| < |P_{n,k}\langle t \rangle|$. 

Now suppose the sets in $\mathcal{A}$ do not have $t$ common elements. So $\mathcal{A}$ is a non-trivial $t$-intersecting family. By Lemma~\ref{lemma3}, there exists a set $J$ such that $|J| \leq 3k-2t-1$ and $|X \cap J| \geq t+1$ for any $X \in \mathcal{A}$. So $\mathcal{A} \subseteq \bigcup_{T \in {J \choose t+1}} \{X \in \mathcal{P}_{n,k} \colon T \subset X\}$. Let $T^* \in {J \choose t+1}$ such that $|\{X \in \mathcal{P}_{n,k} \colon T \subset X\}| \leq |\{X \in \mathcal{P}_{n,k} \colon T^* \subset X\}|$ for all $T \in {J \choose t+1}$. Let $\mathcal{B} = \{X \in \mathcal{P}_{n,k} \colon T^* \subset X\}$. We have
\begin{align} |\mathcal{A}| &\leq \left| \bigcup_{T \in {J \choose t+1}} \{X \in \mathcal{P}_{n,k} \colon T \subset X\} \right| \leq \sum_{T \in {J \choose t+1}} |\{X \in \mathcal{P}_{n,k} \colon T \subset X\}| \leq \sum_{T \in {J \choose t+1}} |\mathcal{B}| \nonumber \\
&\leq {|J| \choose t+1} |\mathcal{B}| \leq {3k-2t-1 \choose t+1} |\mathcal{B}|. \nonumber
\end{align}
%
Let $B = \{{\bf b} \in P_{n,k} \colon f({\bf b}) = X \mbox{ for some } X \in \mathcal{B}\}$. Since $f$ is bijective, $|B| = |\mathcal{B}|$. Let $(e_1,i_1), \dots, (e_{t+1},i_{t+1})$ be the elements of $T^*$. Then, by definition of $B$, $t+1$ of the entries of each member of $B$ are $e_1, \dots, e_{t+1}$. So $|B| \leq p_{n - q, k - (t+1)}$, where $q = e_1 + \dots + e_{t+1}$. Thus, since $q \geq t+1$, $|B| \leq p_{n-(t+1),k-(t+1)}$ by Lemma~\ref{lemma1}.  

Let $n' = n-t$ and $k' = k-t$. 
Since $n' \geq {3k-2t-1 \choose t+1}k^3 - t > {3k-2t-1 \choose t+1}(k')^3$, we have $p_{n',k'} > {3k-2t-1 \choose t+1}p_{n',k'-1}$ by Lemma~\ref{lemma2}. By Lemma~\ref{lemma1}, $p_{n'-1,k'-1} \leq p_{n',k'-1}$. Thus, since $|A| = |\mathcal{A}| \leq {3k-2t-1 \choose t+1}|\mathcal{B}|$ and $|\mathcal{B}| = |B| \leq p_{n'-1, k'-1}$, we have $|A| \leq {3k-2t-1 \choose t+1} p_{n'-1, k'-1} < p_{n',k'}$. Since $|P_{n,k}\langle t \rangle| = p_{n',k'}$, the result follows.~\hfill{$\Box$}\\
\\
\textbf{Proof of Theorem~\ref{main3}.} Let $k \geq t+2$ and $n \geq {3k-2t-1 \choose t+1}k^3$. Let $A$ be a properly $t$-intersecting subset of $P_{n,k}$ such that $A \neq P_{n,k}([t])$. We prove the result by showing that $|A| < |P_{n,k}([t])|$.

For each ${\bf a} = (a_1, \dots, a_k) \in P_{n,k}$, let $R_{\bf a} = \{a_i \colon i \in [k]\}$. Let $\mathcal{A} = \{R_{\bf a} \colon {\bf a} \in A\}$. So $|X| \leq k$ for each $X \in \mathcal{A}$. Since $A$ is a $t$-intersecting set, $\mathcal{A}$ is a $t$-intersecting family.

Suppose the sets in $\mathcal{A}$ have a common $t$-element subset $T$ (i.e.~$T \subseteq \bigcap_{X \in \mathcal{A}} X$). Then $A \subseteq P_{n,k}(T)$. So $|A| \leq |P_{n,k}(T)| = p_{n-s,k-t}$, where $s = \sum_{a \in T} a$. Now $s \geq \sum_{i \in [t]} i = t(t+1)/2$, and equality holds only if $T = [t]$. If $s > t(t+1)/2$, then, by Lemma~\ref{lemma1}, $|A| < p_{n-t(t+1)/2,k-t} = |P_{n,k}([t])|$. Suppose $s = t(t+1)/2$. Then $T = [t]$. So $A \subset P_{n,k}([t])$. Since $A \neq P_{n,k}([t])$, $|A| < |P_{n,k}([t])|$. 

Now suppose the sets in $\mathcal{A}$ do not have a common $t$-element subset. So $\mathcal{A}$ is a non-trivial $t$-intersecting family. By Lemma~\ref{lemma3}, there exists a set $J$ such that $|J| \leq 3k-2t-1$ and $|X \cap J| \geq t+1$ for any $X \in \mathcal{A}$. So $A \subseteq \bigcup_{T \in {J \choose t+1}} P_{n,k}(T)$. Let $T^* \in {J \choose t+1}$ such that $|P_{n,k}(T)| \leq |P_{n,k}(T^*)|$ for all $T \in {J \choose t+1}$. Let $q = \sum_{a \in T^*}a$ and $r = \sum_{i \in [t+1]} i = (t+1)(t+2)/2$. Then $|P_{n,k}(T^*)| = p_{n-q,k-t-1}$, $q \geq r$, and hence, by Lemma~\ref{lemma1}, $|P_{n,k}(T^*)| \leq p_{n-r,k-t-1}$. Therefore, we have
\begin{align} |A| &\leq \left| \bigcup_{T \in {J \choose t+1}} P_{n,k}(T) \right| \leq \sum_{T \in {J \choose t+1}} |P_{n,k}(T)| \leq \sum_{T \in {J \choose t+1}} |P_{n,k}(T^*)| \nonumber \\
&\leq {|J| \choose t+1} p_{n-r,k-t-1} \leq {3k-2t-1 \choose t+1} p_{n-r,k-t-1}. \nonumber
\end{align}
Let $n' = n-t(t+1)/2$ and $k' = k-t$. 
Since $n' \geq {3k-2t-1 \choose t+1}k^3 - t(t+1)/2 > {3k-2t-1 \choose t+1}(k')^3$, we have $p_{n',k'} > {3k-2t-1 \choose t+1}p_{n',k'-1}$ by Lemma~\ref{lemma2}. By Lemma~\ref{lemma1}, $p_{n-r,k-t-1} \leq p_{n',k'-1}$. Thus, since $|A| \leq {3k-2t-1 \choose t+1}p_{n-r,k-t-1}$, we have $|A| \leq {3k-2t-1 \choose t+1}p_{n',k'-1} < p_{n',k'}$. Since $|P_{n,k}([t])| = p_{n',k'}$, the result follows.~\hfill{$\Box$}\\

\end{document}